\documentclass[12pt]{amsart}
\usepackage{amsmath,amstext,amsfonts,amsthm,amssymb}
\usepackage{eucal}
\renewcommand{\subsubsection}[1]{\addtocounter{subsubsection}{1}
{\ \\[3pt]\bf \thesubsubsection. \  #1} }

\swapnumbers

{  \theoremstyle{definition}

}
\newcommand{\hra}{\hookrightarrow}
\newcommand{\lra}{\longrightarrow}

\newcommand{\dpar}{\partial}

\newcommand{\tf}{\tilde f}
\newcommand{\tih}{\tilde h}
\newcommand{\tnabla}{\tilde\nabla}

\newcommand{\CO}{\mathcal{O}}

\newcommand{\fM}{\frak{M}}

\newcommand{\BC}{\mathbb{C}}

\newcommand{\BN}{\mathbb{N}}
\newcommand{\BP}{\mathbb{P}}
\newcommand{\BR}{\mathbb{R}}
\newcommand{\BZ}{\mathbb{Z}}

\begin{document}


\centerline{HILBERT SERIES AND BEREZIN-GELFAND DUALITY}

\vspace{1cm}

\centerline{Vadim Schechtman}

\vspace{2cm}

\centerline{\bf Introduction}

\bigskip\bigskip

{\bf I.1.} In their great work on spherical functions [BG] Berezin and Gelfand wrote: 

{\it "... there exists a deep duality between the function ... giving the law of multiplication in the center of the} [infinitesimal] {\it group ring} [of a semisimple Lee 
group] {\it and the function ... giving multiplication of representations. 

... an analogous duality exists between matrix elements of ... an irreducible representation 
of the group $SU(2)$ ... and the so-called "Clebsch-Gordan coefficients"... 

Another example of such a duality are the formulas of Gelfand-Tsetlin 
for matrix elements of irreducible representations of the algebra of complex 
matrices with trace $0$ and the formulas for coordinates in the group of unitary 
matrices... 
In all of these cases the duality consists in the fact that functions of discrete 
arguments satisfy finite difference equations analogous to differential equations 
satisfied by functions of real variables that correspond to them."} 

The second of the above examples may be expressed by saying that 
we have a duality between the  
classical orthogonal polynomials (Jacobi etc.) and their discrete analogues 
(Hahn etc.). (In fact, all of the above examples admit a similar 
reformulation.)

The main purpose of the present note is to propose an example illustrating that exactly this type of dual 
polynomials appears in certain Hilbert series. Namely, let us say that 
two polynomials $Q(s)$ and $P(t)$ of the same degree $d$ and such that 
$Q(0) = P(0) = 1,\ P(1) \neq 0$, are {\it Euler dual} if there is an equality 
of formal power series  
$$
\sum_{n=0}^\infty\ Q(n)t^n = \frac{P(t)}{(1-t)^{d+1}}
\eqno{(I1)}
$$
(we will see shortly that the equality of degrees is equivalent to 
$Q(-1)\neq 0$). 

Here is an example of an Euler dual pair, which is the main observation of this note. 

{\it Theorem} 1. For each integer $m\geq 0$ we have
$$
\sum_{n=0}^\infty\ Q_m(n)t^n = \frac{P_m(t)}{(1-t)^{m+1}}
\eqno{(I2)}
$$
where 
$$
P_m(t) = (m+1)^{-1}(t-1)^m P_m^{(1,1)}((t+1)/(t-1))
$$
and 
$$
Q_m(s) = [(m+1)!]^{-1}h_m^{(1,1)}(s-1,-2)
$$
Here $P_m^{(\alpha,\beta)}(x)$ denote the {\it Jacobi polynomials} and 
$h_m^{(\alpha,\beta)}(s,N)$ denote the {\it Hahn polynomials} which are the discrete 
analogues of $P_m^{(\alpha,\beta)}(x)$ (their definitions are recalled below).  

We have $\deg\ Q_m = \deg\ P_m = m$; 
$P_m(0) = Q_m(0) = 1$.
 
These polynomials satisfy the following properties: (a)\ 
All coefficients of $P_m(t)$ are positive integers, $P_m(1/t) = t^{-m}P_m(t)$, and 
all its roots are situated on the real half-line $-\infty < t < 0$. 
(b)\ If $n$ is an integer then $Q_m(n)$ is an integer, $Q_m(-1-s) = (-1)^mQ_m(s)$, 
and all its roots  lie on the vertical line 
$\Re\ s = - 1/2$. 
(Formally, the polynomials $Q_m(s)$ may be interpreted as spherical functions 
on a (non-existent) homogeneous space $(\Sigma_{-2}\times \Sigma_{-2})/\Sigma_{-2}$ where $\Sigma_N$ denotes the symmetric group on $N$ letters, 
cf. [D].) 

The geometric meaning of the Hilbert series (I2) is as follows. Consider 
the Grassmanian $X_m = Gr(2,m+3)\subset \BP^{d_m}$ of two-dimensional planes in 
$\BC^{m+3}$ embedded by Pl\"ucker into the projective space, $d_m = (m+3)(m+2)/2$. 
Let $Y_m = X_m\cap L \subset L\cong \BP^{d_{m-1}}$ be its section by a 
generic linear subspace $L\subset \BP^{d_m}$ of codimension $m+2$. 

{\it Theorem} 2. The series (I1) is the Hilbert series of the embedding 
$i_m:\ Y_m\subset L$, i.e. $Q_m(n) = \dim H^0(Y_m,i_m^*\CO_L(n))$, $m,n\geq 0$. 

As a consequence we get another proof of an elegant result due to David Beckwith 
[B] (cf. also [BK]). In fact this remarkable paper was the starting point of the present note.  
Our Theorem 2 is the discrete, or Gelfand dual counterpart 
of Beckwith's theorem. 
As a second consequence, we get another proof of [Bran], Corollary 7.2 
for the series $A_m$. 


The above relation between {\it Euler duality} and {\it Gelfand duality} 
is the first idea of this note. 

{\bf I.2.} The other idea which we wanted to discuss is analogy of Euler duality to {\it Mellin transform}. It can be immediately seen already on the formal level: 

Mellin transform: 
$$
Mf(s) =  \int_0^\infty\ f(x)x^{s-1}dx = \phi(s)
$$
Inverse: 
$$
M^{-1}\phi(x) = \frac{1}{2\pi i}\int_{c-i\infty}^{c-i\infty}\ 
\phi(s)x^{-s}ds = f(s)
$$
We see that $M^{-1}$ resembles a continuous analogue to generating function 
(i.e. to Euler transform), whereas $M$ resembles  taking coefficients of a 
power series,  i.e. to inverse Euler:

\begin{tabular}{ccc}

$\phi(s) = \int_0^\infty\ f(x)x^{s-1}dx$ & \ & 
$h(n) = \frac{1}{2\pi i}\oint F(t)t^{-n-1}dt$ \\

\ & \ &\ \\

$f(x) = \frac{1}{2\pi i}\int_{c-i\infty}^{c-i\infty}\ 
\phi(s)x^{-s}ds$ & \ & $F(t) = \sum_{n=0}^\infty\ h(n)t^n$ \\

\end{tabular} \\

We see similarities and differences. 

Certainly this analogy is behind the scene in [RV]. 
We review in the first section the main features of this analogy (in particular, 
the analogue of Hecke lemma, which in this case is a theorem of 
Popoviciu - Ehrhart - Stanley).

In the last section we introduce some polynomials generalising the above ones. 
They are enumerated by pairs $(\lambda,\epsilon)$ where $\lambda$ is a Young diagram and $\epsilon = \pm$; we discuss their properties and make some conjectures 
about them.

The useful discussions with V.Hinich, V.Gorbounov and F.Hirzebruch 
are gratefully acknowledged.  

This work was finished during the author's stay at Max-Planck-Institut 
f\"ur Mathematik in Bonn, in August 2009.   

\bigskip\bigskip



\centerline{\bf \S 0. Jacobi and Hahn polynomials (recollections)}

\bigskip\bigskip

{\bf 0.1.} {\it Finite differences.} $(a)_n = a(a+1)\cdot (a+n-1) = \Gamma(a+n)/\Gamma(a)$

$\Delta f(x) = f(x+1) - f(x)$; $\nabla f(x) = f(x) - f(x-1)$

$\Delta(fg)(x) = \Delta f(x)\cdot g(x+1) + f(x)\Delta g(x)$

{\bf 0.2.} {\it Hypergeometric function.} 
$$
_pF_q(a_1,\ldots,a_p;b_1,\ldots,b_q;x) = 
\sum_{n=0}^\infty\ \frac{(a_1)_n\ldots (a_p)_n}{(b_1)_n\ldots (b_q)_n}
\frac{x^n}{n!}
$$
{\it Gauss hypergeometric function} $F :=\ _2F_1$

{\bf 0.3.} {\it Jacobi polynomials} (cf. [BE], 10.8, [NU]) $P_n^{(\alpha,\beta)}(x),\ 
n=0,1,\ldots$, are polynomials orthogonal on $[-1,1]$ 
with respect to the scalar product 
$$
(f,g) = \int_{-1}^1\ f(x)g(x)(1-x)^\alpha(1+x)^\beta dx
$$
They satisfy the differential equation
$$
(1-x^2)y'' + (-(\alpha+\beta+2)x+\beta-\alpha)y'' + 
(-n^2+n(\alpha+\beta+1))y = 0
\eqno{(0.3.0)}
$$
They can be defined by the {\it Rodrigues formula}
$$
P_n^{(\alpha,\beta)}(x) = \frac{(-1)^n}{2^n n! (1-x)^\alpha(1+x)^\beta}
\frac{d^n}{dx^n}[(1-x)^{\alpha+n}(1+x)^{\beta+n}]
\eqno{(0.3.1)}
$$   
We have:
$$
P_n^{(\alpha,\beta)}(x) = \frac{1}{2^{n}}\sum_{k=0}^n\ 
\binom{n+\alpha}{k}\binom{n+\beta}{n-k}(x-1)^{n-k}(x+1)^{k}
\eqno{(0.3.2)}
$$
$$
P_n^{(\alpha,\beta)}(x) = \binom{n+\alpha}{n} 
F(-n,n+\alpha+\beta+1;\alpha+1;(1-x)/2)
\eqno{(0.3.3)}
$$
$$
P_n^{(\alpha,\beta)}(-x) = (-1)^nP_n^{(\beta,\alpha)}(x)
\eqno{(0.3.4)}
$$
$$
\frac{d}{dx}P_n^{(\alpha,\beta)}(x) = \frac{n+\alpha+\beta+1}{2}
P_{n-1}^{(\alpha+1,\beta+1)}(x)
\eqno{(0.3.5)}
$$

{\it Particular cases.} {\it Gegenbauer (ultraspherical) polynomials:}
$$
C_n^\lambda(x) = P_n^{(\lambda-1/2,\lambda-1/2)}(x)
$$ 

{\it Legendre polynomials:} $P_n(x) = P_n^{(0,0)}(x) = C_n^{-1/2}(x)$

Generating function:
$$
\sum_{n=0}^\infty\ P_n(x)y^n = (1-2xy+y^2)^{-1/2}
\eqno{(0.3.6)}
$$
Derivation:
$$
(x^2-1)P'_n(x) = n[xP_n(x) - P_{n-1}(x)]
\eqno{(0.3.7)}
$$

\bigskip

{\it Hahn polynomials}

\bigskip

{\bf 0.4.}
Let us recall first some classical definitions and formulas 
related to orthogonal polynomials of a discrete variable, cf.  
[H], [WE],  [NSU].

We consider polynomial solutions $f(x)$ of a finite difference equation 
$$  
\sigma(x)\Delta\nabla f(x) + \tau(x)\Delta f(x) + \lambda f(x) = 0
\eqno{(0.4.1)}
$$
where $\sigma(x), \tau(x)\in\BC[x]$, $\deg\ \sigma(x)\leq 2, 
\deg\ \tau(x)\leq 1$, $\lambda\in\BC$. 

Given $\sigma(x),\ \tau(x)$ as above, let $\rho(x)$ be function (not 
necessarily polynomial!) satisfying  
to the difference equation 
$$
\nabla(\sigma\rho) = \tau\rho
\eqno{(0.4.2)}
$$
Set
$$
\lambda_n = - n\tau' - \frac{n(n-1)}{2}\sigma'',
\eqno{(0.4.3)} 
$$
$n = 0, 1, \ldots$. Then 
$$
f_n(x) = \frac{B_n}{\rho(x)}\nabla^n[\rho(x+n)\prod_{k=1}^n\sigma(x+k)] 
= \frac{B_n}{\rho(x)}\Delta^n[\rho(x)\prod_{k=0}^{n-1}\sigma(x-k)]
\eqno{(0.4.4)}
$$
is a polynomial of degree $\leq n$ satisfying (0.4.1) with 
$\lambda = \lambda_n$. The number $B_n$ is a normalising constant: 
if we want that $f_n(x) = x^n + O(x^{n-1})$ then 
$$
B_n = \prod_{k=0}^{n-1}\bigl(\tau' + \frac{n+k+1}{2}\sigma'')^{-1}
$$
(so if we want the existence of a solution of degree $n$, 
the factors in the product should be $\neq 0$; recall that $\deg\tau\leq 1$ and $\deg\sigma\leq 2$, so $\tau'$ and $\sigma''$ are constants). Formula $(0.4.4)$ is the finite 
difference analogue of the Rodrigues formula. 

{\bf 0.5.} More specifically, let $\alpha, \beta, N$ be $3$ numbers. 
Set 
$$
\sigma(x,\alpha,\beta,N) = x(N+\alpha-x)
\eqno{(0.5.1a)}
$$
$$
\tau(x,\alpha,\beta,N) = -(2+\alpha+\beta)x + (N-1)(\beta+1) 
\eqno{(0.5.1b)}
$$  
$$
\rho(x,\alpha,\beta,N) = \frac{\Gamma(N+\alpha-x)\Gamma(x+\beta+1)}{\Gamma(x+1)\Gamma(N-x)}
\eqno{(0.5.1c)}
$$
and $B_n = (-1)^n/n!$. 
Then the polynomials given by (0.4.4) are called {\it the Hahn polynomials}\footnote{Wolfgang Hahn (1911 - 1998), a student of Issai Schur, PhD 1933.} 
and denoted by $h_n^{(\alpha,\beta)}(x,N)$. 

Explicitly, 
$$
h_m^{(\alpha,\beta)}(n,N) = 
$$
$$
\frac{(-1)^m(N-1)!(\beta+1)_m}
{m!(N-m-1)!}\ _3F_2(-m,\alpha+\beta+m+1,-n;\beta+1,1-N;1)
\eqno{(0.5.2)}
$$
cf. [NSU], (75a). The first polynomials are: 
$$
h_0^{(\alpha,\beta)}(n,N) = 1;\ h_1^{(\alpha,\beta)}(n,N) = 
(2+\alpha+\beta)x - (N-1)(\beta+1) 
\eqno{(0.5.2b)}
$$
cf. (0.5.1b).

If $N$ is a positive integer then the $h_m$'s satisfy the othogonality relation
$$
\sum_{n=0}^{N-1}\ h_m^{(\alpha,\beta)}(n,N)h_{m'}^{(\alpha,\beta)}(n,N)
\rho(n,\alpha,\beta,N) = d^2_m(\alpha,\beta,N)\delta_{m,m'}
\eqno{(0.5.3)}
$$
for suitable constants $d^2_m(\alpha,\beta,N)$. 

{\it Relation to Jacobi polynomials}:
$$
\lim_{N\rightarrow\infty}\ N^{-m}h_m^{(\alpha,\beta)}(Nx,N) = 
P_m^{(\alpha,\beta)}(2x-1)
\eqno{(0.5.4)}
$$
cf. [NSU], (56).

{\bf 0.6.}
We have 
$$
\Delta h_m^{(\alpha,\beta)}(n,N) = (\alpha+\beta+m+1)h_{m-1}^{(\alpha+1,\beta+1)}(n,N-1)
\eqno{(0.6.1)}
$$
$$
h_m^{(\alpha,\beta)}(N-1-x,N) = (-1)^mh_m^{(\beta,\alpha)}(x,N)
\eqno{(0.6.2)}
$$

\bigskip\bigskip



\centerline{\bf \S 1. Euler transform and toy Hecke lemma}

\bigskip\bigskip

{\bf 1.1.} {\it Euler transform.} $Q(s)\in \BC[s]$, $\deg Q = d$, $Q(0) = 1$ whence 
$$
\sum_{n=0}^\infty\ Q(n)t^n = \frac{P(t)}{(1-t)^{d+1}}
$$
where $P(t) \in \BC[t]$, $\deg P = e \leq d,\ P(1)\neq 0,\ 
P(0) = 1$. 
Let us say that $P(t)$ is the {\it Euler transform} of $Q(s)$, and write $P(t) = EQ(t)$.  

We will see that this operation is in many respects  analogous to {\it inverse Mellin transform}. 

Let us call $f(Q) = d-e$ the {\it defect} of $Q$. 

For a natural $a$ denote $h_a(s) = (s+1)\ldots (s+a)$. 

{\bf 1.1.1.} $f(Q) = \text{max}\{a|\ h_a(s)|Q(s)\}$. 
In other words, $f(Q) = \text{max}\{a|\ Q(-1) = Q(-2) = \ldots = Q(-a) = 0$. 

We define $R$ by
$$
\sum_{n=0}^\infty\ R(n)t^n = \frac{P(t)}{(1-t)^{e+1}}  
$$
Then $\deg R = e$, $R(0) = 1$. We say that $R$ is the {\it inverse Euler transform} 
of $P$ and write $R(s) = \tilde EP(s)$.  

$$
(1-t)\sum_{n=0}^\infty\ S(n)t^n = \sum_{n=0}^\infty\ \tilde\nabla S(n)t^n
$$
where $\tilde\nabla S(n) = S(n) - S(n-1)$ if $n>0$ and $S(0)$ for $n=0$. 
If $S(-1)=0$ then $\tnabla S(n) = \nabla S(n)$ for $n\geq 0$. 

Here $\nabla f(x) = f(x) - f(x-1)$.  

It follows: 

{\bf 1.1.2.} $R = \nabla^{f(Q)}Q$

\bigskip

{\bf 1.2.} Set 
$$
F(t) = \sum_{n=0}^\infty\ Q(n)t^n = \frac{P(t)}{(1-t)^{d+1}}
$$
Let us consider $F(t)$ as a rational function on $\BC$, so $F(t^{-1})$ 
is also a rational function, and we can consider its Taylor expansion at $0$.
 
("Taylor", not Laurent, since $F(t)$ is regular at infinity, moreover, it has $0$ of order $f(Q)+1$ at $\infty$.) 

{\bf 1.2.1.} {\it Theorem} (Tiberiu Popoviciu, cf. [S], 4.6). $F(t^{-1}) = 
\sum_{n=1}^\infty\ Q(-n)t^n$. 

We see again that $f(Q) + 1 = \text{ord}_0 F(t^{-1}) = \text{ord}_\infty F(t)$.  

{\it Proof.} First verify this for 
$$
F(t) = \frac{1}{(1-t)^{p+1}} = \sum_{n\geq 0}\binom{n+p}{p}t^n
$$
The general case follows by Taylor expansion of $Q$ at $t=1$.

{\bf 1.3.} {\it Another interpretation.} For a polynomial $f(s)$ we have 
its discrete Taylor expansion at $x=-1$:
$$
f(s) = \sum_{a=0}^\infty\ \frac{\nabla^af(-1)}{a!}h_a(s)
$$
So we see that 

$\text{max}\{a:\ h_a(s)|f(s)\} = \text{min}\{b:\ \nabla^bf(-1)\neq 0\}$; 
one can call this number the "`discrete order of zero" of $f$ at $s=-1$; 
let us denote it $\text{ordd}_{s=-1}f(s)$.

Theorem of Popoviciu implies that the 
$$
\text{ordd}_{s=-1}Q(s) = \text{ord}_\infty F(t)
$$  
  
{\bf 1.4.} {\it Corollary} ("toy Hecke lemma", cf. [S], 4.7). Let $\epsilon\pm 1$, $f=f(Q)$. Then 
$$
P(t^{-1}) = \epsilon t^{-e}P(t)
$$
iff 
$$
Q(m) = \epsilon (-1)^d Q(-f-1-m)
$$
for all $m\in \BZ$.  

(Note that "in the critical strip" $Q(-1) = \ldots = Q(-f) = 0$ by our hypothesis.) 

{\it Proof.} $P(t^{-1}) = \epsilon t^{-e}P(t)$ iff 
$F(t^{-1}) = (-1)^{d+1}\epsilon t^{d+1-e}F(t) = \epsilon (-1)^{d+1} t^{f+1}F(t)$. By Theorem, 
$$
F(t^{-1}) = - \sum_{k=1}^\infty Q(-k)t^k = 
\epsilon (-1)^{d+1} \sum_{n=0}^\infty Q(n)t^{n+f+1} = 
\epsilon (-1)^{d+1} \sum_{k=f+1}^\infty Q(k-f-1)t^k
$$
so we get our assertion by putting $m=-k$.

\bigskip\bigskip



\centerline{\bf \S 2. Three sets of polynomials} 

\bigskip\bigskip

{\bf 2.2.} Let us define the following sequences of polynomials, indexed by natural
numbers:  

(a)  
$$
f_m(t) = \frac{1}{m+2} \sum_{k=0}^m \binom{m}{m-k}\binom{2m+2-k}{m+1-k} t^k = 
\sum_{k=0}^m\ a_{km} t^k
$$
For example: 
$$
f_0 = 1,\ f_1 = 2 + t,\ f_2 = 5 + 5t + t^2,\ f_3 = 14 + 21 t + 9 t^2 +
t^2,\ldots 
$$
We set
$$
\tf_m(t) = t^mf_m(t^{-1}) = 
\frac{1}{m+2} \sum_{k=0}^m \binom{m}{k}\binom{m+2+k}{1+k} t^k  
$$
These polynomials have been known since long ago, cf. [K]. 

(b)
$$
g_m(t) = \frac{1}{m+1} \sum_{j=0}^m \binom{m+1}{j+1}\binom{m+1}{j} t^j 
= \sum_{k=0}^m\ b_{km} t^k
$$
For example: 
$$
g_0 = 1,\ g_1 = 1 + t,\ g_2 = 1 + 3t + t^2,\ g_3 = 1 + 6t + 6t^2 + t^3, 
\ldots
$$

{\bf 2.2.} {\it Theorem} (a) 
$$
\tf_m(t) = F(-m,m+3;2;-t) = \frac{1}{m+1}P_m^{(1,1)}(2t+1)
\eqno{(2.2.1)} 
$$
$$
f_m(t) = \frac{t^m}{m+1}P_m^{(1,1)}\bigl(\frac{t+2}{t})
\eqno{(2.2.2)}
$$

(b) 
$$
g_m(t) = \frac{(t-1)^m}{m+1}P_m^{(1,1)}\bigl(\frac{t+1}{t-1}\bigr)
\eqno{(2.2.3)} 
$$

{\it Proof.} (a) The first equality in (2.2.1) follows from the definition; the second one 
--- from (0.3.3). The equality (2.2.2) immediately follows.  

(b) Take (0.3.2) with $x = (t+1)/(t-1)$ and $\alpha = \beta = 1$. 

{\bf 2.2.2.} {\it Corollary} (cf. [B], formula (6) and Theorem). Let $P_{m+1}(x)$ 
denote the Legendre polynomial. Then  
$$\tf_m(t) = \frac{2}{(m+1)(m+2)}P'_{m+1}(2t+1)$$

This is equivalent to (2.2.1): it suffices to take into account that
$$
P'_{m+1}(x) = \frac{m+2}{2}P_m^{(1,1)}(x)
$$
by (0.3.5).

Just for completeness we add 

{\bf 2.2.2.} {\it Corollary} (D.Beckwith, [B]). The generating function
$$
f(x,y) = \sum_{n=3}^\infty\ \tf_{n-3}(x)y^{n-1}
\eqno{(2.2.4)}
$$
satisfies the differential equation with initial condition:  
$$
f'_x = f f'_y,\ f(0,y) = y^2(1-y)^{-1}
\eqno{(2.2.5)}
$$

{\it Proof (op. cit.)} One verifies directly that if the function $\phi(x,y)$ 
satisfies the functional equation 
$$
\phi(1-y-x\phi) = (y+x\phi)^2
\eqno{(2.2.6)}
$$
then it satisfies (2.2.5). Explicitly, the solution to (2.2.6) is  
$$
\phi(x,y) = \frac{2(1-\rho y - \sqrt{1-2\rho y + y^2})}{\rho^2-1}
$$
where $\rho = 2x+1$. Now, using (0.3.6) and (0.3.7) one sees that 
$$
\phi(x,y) = \sum_{n=3}^\infty\ \frac{2P'_{n-2}}{(n-1)(n-2)}y^{n-1}
$$
and now 2.2.1 implies that $\phi(x,y) = f(x,y)$.

{\bf 2.3.} {\it Corollary.} The polynomials $g_m(t)$ are {\it self-reciprocal}, 
$g_m(t^{-1}) = t^{-m}g_m(t)$. 

This follows from (2.2.3) and (0.3.4). Of course one can see this immediately 
from the definition of $g_m(t)$. 

{\bf 2.4.} {\it Corollary} (cf. [L], Theorem 3) $g_m(t) = f_m(t-1)$

{\it Proof.} Replace $t$ by $t-1$ in (2.2.2) and you get (2.2.3). This also may be 
verified directly. 

\bigskip

{\bf 2.5.} We define another sequence of polynomials:
$$
h_m(n) := \frac{(n+1)(n+2)^2\ldots (n+m+1)^2(n+m+2)}
{1\cdot 2^2\ldots (m+1)^2\cdot (m+2)},
\eqno{(2.5.1)}
$$
$m= 0, 1, 2, \ldots$. Note that $\deg\ h_m = 2m+2$. 

The following theorem is well known. 

{\bf 2.5.2.} {\it Theorem.} 
$$
\sum_{n=0}^\infty\ h_m(n)t^n = \frac{g_m(t)}{(1-t)^{2m+3}}
\eqno{(2.5.2)}
$$

{\it Proof.} The following elementary lemma is usefull (cf. [GW]): 

{\bf 2.5.2.} {\it Lemma.} If $\phi(x)\in\BC[x]$ then 
$$
\sum_{n=0}^\infty\ \phi(n)t^n = \phi(t\partial_t)\biggl(\frac{1}{1-t}\biggr)
$$

{\bf 2.5.3.} {\it Lemma.} (a) 
$$
\biggl(1 + t\dpar_t\biggr)\biggl(1 + \frac{t\dpar_t}{2}\biggr)\ldots 
\biggl(1 + \frac{t\dpar_t}{\ell}\biggr)\frac{1}{1-t} = 
\frac{1}{(1-t)^{\ell+1}}
$$
(b) 
$$
\biggl(1 + \frac{t\dpar_t}{2}\biggr)\ldots 
\biggl(1 + \frac{t\dpar_t}{i+1}\biggr)\frac{1}{(1-t)^{\ell+1}} = 
\frac{\sum_{j=0}^i \frac{1}{\ell}\binom{i}{j}\binom{\ell}{j+1}t^j}
{(1-t)^{\ell+i+1}}
$$

{\it Proof.} Induction on $\ell$, then on $i$. 

To finish the proof of 2.6, we remark that $h_m = h_m^{(2)}h_m^{(1)}$ where 
$$
h_m^{(1)}(n) := \frac{(n+1)(n+2)\ldots (n+m+1)(n+m+2)}
{(m+2)!}
$$
and
$$
h_m^{(2)}(n) := \frac{(n+2)\ldots (n+m+1)}
{2\ldots (m+1)}
$$
and use 2.5.3 with $\ell = m+2,\ i = m+1$.

\bigskip\bigskip

\newpage

\centerline{\bf \S 3. Hahn polynomials vs Hilbert polynomials} 

\bigskip\bigskip 

{\bf 3.1.} 
Note the following evident properties of the polynomials 
$h_m$: 
$$
h_m(n) = 0\ \text{for\ }n=-1, -2,\ldots -m-2
\eqno{(3.1.1)}
$$
$$
h_m(-m-3-n) = h_m(n)
$$

Let us 
define polynomials $Q_m$ by the generating series 
$$
\sum_{n=0}^\infty\ Q_m(n)t^n = \frac{g_m(t)}{(1-t)^{m+1}}
\eqno{(3.1.2)}
$$
So we have $\deg Q_m = m$,  
$$
Q_m(n) = \nabla^{m+2}h_m(n)
\eqno{(3.1.3)}
$$
Note that if $f(a-x) = cf(x)$ then $\nabla(a+1-x) = - c\nabla f(x)$. 
It follows that 
$$
Q_m(-1-n) = (-1)^mQ_m(n)
\eqno{(3.1.4)}
$$

{\bf 3.2.} {\it Theorem.} 
$$
Q_m(n) = \frac{1}{(m+1)!}h_m^{(1,1)}(n-1,-2)
\eqno{(3.2.1)}
$$

{\it Proof.} Let us set in the discussion 0.5 $\alpha = \beta = 0$, 
$N = -1$. Then we obtain: 
$$
\sigma(x) = - x(x+1),\ \rho(x) = 1,\ 
\tau(x) = \Delta\sigma(x) = -2x
\eqno{(3.2.2)}
$$
It follows that
$$
h_{m+1}^{(0,0)}(n,-1) = \frac{\nabla^{m+1}[
(n+1)(n+2)^2\ldots (n+m+1)^2(n+m+2)]}{(m+1)!}
$$
whence
$$
Q_m(n) = \frac{1}{(m+2)!}\nabla h_{m+1}^{(0,0)}(n,-1) = 
\frac{1}{(m+2)!}\Delta h_{m+1}^{(0,0)}(n-1,-1) =
$$
$$
\frac{1}{(m+1)!}h_m^{(1,1)}(n-1,-2)
$$
by (3.3.1).  

This theorem should be compared with 2.2: note that $h_m^{(1,1)}$ 
are discrete analogues of $P_m^{(1,1)}$. In fact, the polynomials (3.2.1) 
where introduced (up to a constant multiple) already by Chebyshev, cf. [Ch]. 

\bigskip

\newpage

{\it Geometric meaning of $h_m$'s}

\bigskip

{\bf 3.3.} The result below is a particular 
case of an old theorem due to Hirzebruch.  
Let $Gr(2,m+3)$ be the grassmanian of two-dimensional planes in 
$\BC^{m+3}$. 
Consider the Pl\"ucker embedding
$$
\iota_m:\ Gr(2,m+3) \hra \BP(\Lambda^2\BC^{m+3})\cong \BP^{d_m}
\eqno{(3.3.1)}
$$
where $d_m := (m+3)(m+2)/2$. 
     
Consider the coordinate algebra of $\iota_m$: $A_m = \oplus_{n=0}^\infty\ 
A_m^n$ where 
$$
A_m^n = H^0(Gr(2,m+3),\iota_m^*\CO_{\BP^{d_m}}(n))
$$
Let $H(A_m;t) = \sum_{n=0}^\infty\ \dim\ A_m^n\cdot t^n$ 
be its Hilbert series.  

{\bf 3.4.} {\it Theorem}, cf. [Hir]. $\dim\ A_m^n = h_m(n)$

{\bf 3.5.}  
{\it Proof.} We follow [GW]. Set $\ell = m+2$. We have $Gr(2,\ell+1) = G/P$ where $G = SL(\ell+1)$ and 
$P$ is the obvious parabolic. Let $\lambda = \varpi_2$ be the highest 
weight of the irreducible $G$-module $L(\lambda) = \Lambda^2\BC^{\ell+1}$. 
We can identify $A_m^n$ with the irreducible $G$-module $L(n\lambda)$ 
of highest weight $n\lambda$. 

Its dimension may be calculated using the Hermann Weyl character formula:
$$
\dim\ A_m^n = \dim L(n\lambda) = \prod_{\alpha > 0} \biggl(1 + 
\frac{(\lambda|\alpha)}{(\rho|\alpha)}\cdot n\biggr)
$$
the product over the positive roots of the system $A_\ell$. 
In the notations of Bourbaki, [Bou], the positive roots are 
$\alpha_{ij} = \epsilon_i -
\epsilon_j,\ 1\leq i < j \leq \ell + 1$; the half-sum of the positive roots 
$$
\rho = \frac{1}{2}\sum_{i=0}^\ell\ (\ell - 2i)\epsilon_i
$$
and
$$
\varpi_2 = \epsilon_1 + \epsilon_2 - \frac{2}{\ell + 1}
\sum_{1\leq j\leq \ell+1}\ \epsilon_j
$$
It follows easily that 
$$
\dim L(n\lambda) = \prod_{j=1}^{\ell-1}\bigl(1 + \frac{n}{j}\bigr)\cdot 
\prod_{j=2}^{\ell}\bigl(1 + \frac{n}{j}\bigr) = h_m(n)
$$

{\bf 3.6.} As a corollary we conclude that all numbers $h_m(n)$ are integers for 
$n\in\BN$, hence for all $n\in\BZ$ by (3.1.1). 

Note that the coefficients of $g_m(t)=\sum_{j=0}^m b_{jm}t^j$ are by definition  
$$
b_{jm} = \frac{1}{m+1}\binom{m+1}{j+1}\binom{m+1}{j} = h_{j+1}(m-j) 
$$
(we see here a sort of reciprocity between $j$ and $m$). In particular 
they are all integers.

\bigskip\bigskip



\centerline{\bf 4. Hill polynomials}

\bigskip\bigskip

{\bf 4.1.} A {\it hill} is a sequence of positive integers 
$\mu = (\mu_1,\ldots,\mu_\ell)$ such that 

$\mu_i = \mu_{\ell+1-i}$ and 
$\mu_{i-1}\leq \mu_i$ for $i\leq (\ell+1)/2$. 

$\ell =: w(\mu)$ is called the {\it width} of $\mu$ and $[(\ell+1)/2] =: h(\mu)$ is called the {\it height} of $\mu$.

Another way of looking at hills is as follows. Let $\lambda = (\lambda_1,\ldots, 
\lambda_m),  0<\lambda_1\leq \ldots \leq \lambda_m,\ \lambda_i\in\BZ$, be a Young diagram. Define two hills 
$$
D_+\mu = (\lambda_1,\ldots, \lambda_m,\lambda_m,\ldots,\lambda_1),\ 
D_-\mu = (\lambda_1,\ldots, \lambda_{m-1},\lambda_m,\lambda_{m-1},\ldots,\lambda_1)
$$ 
(the "even" and "odd" doubles of $\lambda$). This way one gets a bijection between 
the set of hills and the set of pairs $(\lambda,\epsilon)$ where 
$\lambda$ is a Young diagram and $\epsilon = \pm$.

Define polynomials 
$$
h_\mu(s) = \prod_{i=1}^\ell\ (s+i)^{\mu_i},\ 
\tih_\mu(s) = h_\mu(s)/h_\mu(0)
$$
and
$$
Q_\mu(s) = \nabla^\ell\tih_\mu(s)
$$
We have $Q_\mu(-1-s) = (-1)^\ell Q_\mu(s)$. These polynomials will be called 
{\it Hill polynomials}. 

{\bf 4.2.} The following discussion is inspired by [RV] (cf. also [G], [M]). For an integer $a\geq 1$ let $\fM_a$ denote the set of hills of width $a$; set $\fM_0 = \{\emptyset\}$.  

Define an operator 
$\dpar:\ \fM_a \lra \fM_{a-1},\ a \geq 1$, by 
$$
(\dpar\mu)_i = \text{min}(\mu_i,\mu_{i+1}),\ 1\leq i \leq a-1
$$
if $a>1$, and $\dpar:\ \fM_1\lra\fM_0$ to be the unique map. Define the sets of corresponding polynomials 
$\Phi_a = \{h_\mu|\mu\in \fM_a\}$, $\Phi_0 = \{1\}$ and an operator 
$\dpar:\ \Phi_a \lra \Phi_{a-1}$ by
$$
\dpar h_\mu = h_{\dpar\mu} = \text{gcd}(h_\mu(x),h_\mu(x-1))
$$
For an arbitrary $\phi\in \fM_a$ 
we denote $\phi_+(x) := \phi(x)/\dpar\phi(x)$ and 
\newline $\phi_-(x) := \phi(x-1)/\dpar\phi(x)$. We have 
$$
\phi_-(x) = \phi_+(a+1-x)
\eqno{(4.2.1)}
$$ 
For any $a\in\BZ$ let $R_a$ denote the set of polynomials  
$p(x)$  all whose roots lie on the line $\Re x = - (a+1)/2$. 

The following Lemma generalises slightly the Lemma from [RV], no. 2. 

{\bf 4.3.} {\it Lemma.} If $\phi\in \Phi_a$ and $p\in R_a$ then 
$\nabla(\phi\cdot p) = \dpar\phi\cdot q$ where $q\in R_{a+1}$. 

{\it Proof} goes along the same lines as in {\it loc. cit.} 
We have $\nabla(\phi p) = (\dpar\phi) q$ where 
$$
q(x) = \phi_+(x) p(x) - \phi_-(x)p(x-1) 
$$ 
If $q(\beta) = 0$ then 
$$
\phi_+(\beta)p(\beta) = \phi_-(\beta)p(\beta-1)
\eqno{(*)} 
$$
Let $p(x) = c'\prod_j\ (x-\mu_j)$, $\phi_+(x) = c\prod_k(x-\nu_k)$, so 
$\phi_-(x) = c\prod_k(-x+a+1-\nu_k)$, cf. (4.2.1).  
For all $j$, $\Re\mu_j = -(a+1)/2$ and $\mu_j\mapsto \mu_j + 1$ establishes 
a bijection between the roots of $p(x)$ (counted with their multiplicities) 
and the roots of $p(x-1)$. Similarly, for all $k$ $\nu_k < - a/2$ and 
$\nu_k \mapsto a+1 - \nu_k$ establishes a bijection between 
the roots of $\phi_+(x)$ and the roots of $\phi_+(x)$. 

In other words, all roots of $\phi_+(x)p(x)$ are situated 
in the right half-plane $\{|z| < - a/2$ and the reflection with respect 
to the line $\Re z = - a/2$ maps them bijectively to the roots of 
of $\phi_-(x)p(x-1)$. It follows that if 
$\Re\gamma < - a/2$ (resp. $> - a/2$) 
then $|\phi_+(\gamma)p(\gamma)|$ is less than (resp. greater than) 
$|\phi_-(\gamma)p(\gamma-1)|$. Thus (*) implies $\Re\beta = - a/2$.

As an immediate corollary we get 

{\bf 4.4.} {\it Theorem}. All roots of $Q_\mu(s)$ lie on the line 
$\Re(s) = - 1/2$. 

{\bf 4.5.} It is likely that $Q_\mu(s)$ satisfies a difference equation of 
order $h = h(\mu)$ of the form 
$$
(\sum_{i=0}^h\ p_i(s)\nabla^i)Q_\mu(s) = 0
$$
where $p_i(s)$ is a polynomial of degree $i$. 

{\it Exercise.} Prove this for hills of height $1$. 

{\bf 4.6.} Define {\it dual hill polynomials} $P_\mu(t)$ by 
$$
\sum_{n=0}^\infty\ Q_\mu(t) t^n = \frac{P_\mu(t)}{(1-t)^{d_\mu+1}}
$$
where $d_\mu = \deg Q_\mu = \deg h_\mu - \ell = \sum_i \mu_i - \ell = \deg P_\mu$. 

We also have
$$
\sum_{n=0}^\infty\ \tih_\mu(t) t^n = \frac{P_\mu(t)}{(1-t)^{v_\mu+1}} 
$$
where $v_\mu := \sum_i \mu_i$. 

We have 
$$
P_\mu(t^{-1}) = t^{-d_\mu}P_\mu(t)
$$

{\bf 4.7.} {\it Conjecture.} All roots of $P_\mu(t)$ are simple and lie on the real 
half line $-\infty < t < 0$.

{\bf 4.8.} {\it Example.} Let $\mu_m$ denote the following hill of width $m+2$ and height $2$: $\mu_m = (1,2,2,\ldots,2,1)$; then evidently 
$Q_{\mu_m} = Q_m$ from (3.1.3), so we conclude that all roots of $Q_m(s)$ lie on the line $\Re s = -1/2$. 

As a corollary of this and of Thm. 3.2, all roots of the Hahn polynomials 
$h_m^{(1,1)}(s,-2)$ lie on the line $\Re s = 1/2$. 

The dual hill polynomials $P_{\mu_m}$ coincide with $g_m$ (note that the coefficients of $f_m(t) = g_m(t+1)$ may be interpreted as numbers of some standard 
Yong tableaux, [S2]).   
Their roots are simple and belong 
to $\BR_{<0}$ by Thm. 2.2 since the roots of Jacobi polynomials are simple and lie in the interval $-1 < t < 1$. 

{\bf 4.9.} {\it Example: Eulerian polynomials.} A sort of "opposite" example is provided by the hills 
of width one: $\nu_k = (k+1)$, $k\geq 0$. In this case we have 
$h_{\nu_k}(s) = \tih_{\nu_k}(s) = (s+1)^{k+1}$, $Q_{\nu_k}(s) = 
(s+1)^{k+1} - s^{k+1}$. Let us denote for brevity $Q_{\nu_k}$ by 
$_kQ$. The roots of $_kQ(s)$ are: 
$$
s_k = - \frac{1}{2} - \frac{i}{2}\cot\biggl(\frac{\pi\ell}{k+1}\biggr),\ \ell=1,\ldots,k
$$
After change of variable $s = -1/2 + ir$ the polynomials $_kQ$ will 
more or less coincide with the polynomials introduced by Euler in his proof 
of the product formula for $\sin x$ and $\cos x$, cf. [Eu] (a), [W], III, \S XIX. 

Passing to the duals 
$$
(1-t)\sum_{n=0}^\infty\ (n+1)^{k+1}t^n = \sum_{n=0}^\infty\ _kQ(n)t^n = 
\frac{P_{\nu_k}(t)}{(1-t)^{k+1}}
$$
we see that $P_{\nu_k}(t)$ conicides with the Eulerian polynomial 
$P_{k+1}(t)$ discussed by Hirzebruch in [Hir2], cf. [Eu] (b). Conjecture 4.7 amounts to saying that all roots of $P_k(t)$ are real and simple. This is a well known fact.

\bigskip\bigskip

\bigskip\bigskip

\newpage

\centerline{\bf References}

\bigskip\bigskip

[BE] Higher transcendental functions, H.Bateman Manuscript Project, 
A.Erd\`elyi (Ed.), McGraw-Hill, 1953

[B] D.Beckwith, Legendre polynomials and polygon dissections?, 
{\it The American Mathematical Monthly} {\bf 105} (1998), 256 - 257.

[BG] F.Berezin, I.M.Gelfand, Some remarks on the theory of spherical functions 
on symmetric Riemannian manifolds, {\it Tr. Mosk. Mat. Ob-va}, 
{\bf 5} (1956), 311 - 351.   

[Bou] N.Bourbaki, Groupes et alg\`ebres de Lie, Chap. VI. Hermann, Paris, 1968. 

[Bran] P.Br\"and\'en, On linear transformations preserving the P\'olya frequency 
property, arXiv:math/0403364. 


[BK] V.M.Buchstaber, E.V.Koritskaya, Quasilinear Burgers-Hopf equation 
and Stasheff polytopes, {\it Funct. Anal. Appl.}, {\bf 41} (2007), 34 - 47. 

[Ch] P.L.Chebyshev, Sur une nouvelle s\'erie, 
{\it Bull. Phys.-Math. de l'Acad. Imp. des Sciences de St. P\'etersbourg}, 
XVII, 257 - 261 = Oeuvres, tome I, 381 - 384. 


[D] Ch.Dunkl, An addition theorem for Hahn polynomials: the spherical functions, 
{\it SIAM J. Math. Anal.} {\bf 9} (1978), 627 - 637. 

[Eu] L.Euler, (a) De summis serierum reciprocarum ex potestatibus numerorum 
naturalium ortarum Dissertatio altera in qua eadem summationes ex fonte 
maxime diverso derivantur, {\it Opera Omnia}, Ser. I, XIV, 138 - 155. 
(b) Remarques sur un beau rapport entre les s\'eries des puissances tant direct 
que r\'eciproques, {\it Ibid.}, XV, 70 - 90.    

[Gol] V.Golyshev, The canonical strip. I, math.AG/0903.2076.  

[GW] B.H.Gross, N.R.Wallach, On the Hilbert polynomials and Hilbert 
series of homogeneous projective varieties, Preprint.  

[H] W.Hahn,\"Uber Orthogonalpolynome, die $q$-Differenzengleichungen 
gen\"ugen, {\it Math. Nachr.} {\bf 2} (1949), 4 - 34.

[Hir] F.Hirzebruch, Characteristic numbers of homogeneous domains, {\it 
Semin. Analytic Functions} {\bf 2} (1958), 92 - 104. 

[Hir2] F.Hirzebruch, Eulerian polynomials, {\it M\"unster J. Math.} {\bf 1} 
(2008), 9 - 14. 

[K] T.P.Kirkman, On the $k$-partitions of the $r$-gon and $r$-ace, {\it Philos. Trans. 
Roy. Soc. London} {\bf 147} (1857), 217 - 272.     

[L] C. W. Lee, The associahedron and triangulations of the $n$-gon, 
{\it Europ. J. Combinatorics} {\bf 10} (1989), 551 - 560.

[Lev] R.J.Levitt, The zeros of Hahn polynomials, {\it SIAM Review} 
{\bf 9} (1967), 191 - 203.

[M] L.Manivel, The canonical strip phenomenon for complete intersections 
in homogeneous spaces, arXiv:0904.2470.   

[NSU] A.F.Nikiforov, S.K.Suslov, V.B.Uvarov, Classical orthogonal 
polynomials of the discrete variable (Russian), Moscow, Nauka, 1985. 

[NU] A.F.Nikiforov, V.B.Uvarov, Special functions of mathematical physics, 
Birkh\"auser, 1988. 


[RV] F.Rodrigues-Villegas, On the zeros of certain polynomials, 
{\it Proc. Amer. Math. Soc.} {\bf 130} (2002), 2251 - 2254. 

[S] R.Stanley, Hilbert functions of graded algebras, {\it Adv. Math.} 
{\bf 28} (1978), 57 - 83.

[S2] R.Stanley, Polygon dissections and standard Young tableaux, 
{\it J. Combinatorial Theory}, Series A {\bf 76} (1996), 175 - 177.  

[WE] Maria Weber, A.Erd\`elyi, On the finite difference analogue 
of Rodrigues' formula, {\it The American Mathematical Monthly} {\bf 59} (1952), 
163 - 168.

[W] A.Weil, Number Theory. An approach through history, Birkh\"auser, 1983.

\end{document}